\documentclass[24pt]{article}

\usepackage{arxiv}

\usepackage[utf8]{inputenc} 
\usepackage[T1]{fontenc}    
\usepackage{hyperref}       
\usepackage{url}            
\usepackage{booktabs}       
\usepackage{amsfonts}       
\usepackage{nicefrac}       
\usepackage{microtype}      
\usepackage{lipsum}         
\usepackage{graphicx}
\usepackage{natbib}
\usepackage{doi}

\usepackage{tikz} 
\usepackage{caption}
\usepackage{amsmath}
\usepackage{cleveref}       

\title{Combinatorial identities using Bernoulli graphs}

\author{ \href{https://orcid.org/0000-0002-8749-3324}{\includegraphics[scale=0.06]{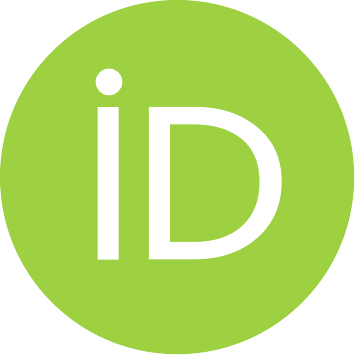}\hspace{1mm}Jacques Bourg}\ }

\renewcommand{\headeright}{}
\renewcommand{\undertitle}{}
\renewcommand{\shorttitle}{Combinatorial identities using Bernoulli graphs}

\hypersetup{
pdftitle={Combinatorial identities using Bernoulli graphs},
pdfsubject={math.CO},
pdfauthor={Jacques Bourg},
pdfkeywords={Combinatorial identities},
}

\begin{document}

\maketitle

\begin{abstract}

In here, I present a series of combinatorial equalities derived using a graph 
based approach. Different nodes in the  graphs are visited following probabilistic dynamics of a moving dot.  The results are presented in such a way that the generalisation (to more nodes, or dimensions) is straightforward.
At an instant m, we 'take a picture' of the system and we compute the probabilities of being at particular positions in space. The sum of all these probabilities is equal to one. 

\end{abstract}


\section{Introduction}
 There are hundreds of combinatorial identities  \citep{gould1972combinatorial}, \citep{Gould2010}, the most known one being the binomial identity.
In here we present some identities obtained using the following method: consider a particle displacing itself given simple   probabilistic dynamics. On certain node configurations, there are multiple ways to attain a point at a given time point. The sum of the probabitilies over space at a given time point must be equal to one.

\section{Chain}

Lets consider a point that advances from the state $c_i$ to the state $c_{i+1}$ with probability $p$, and with a waiting time of one unit. This point reaches also 
the state $c_{i+1}$ with probability $1-p$, but instantaneously. After N states, the point attains a terminal state. 

$$  $$ 

\begin{tikzpicture}[node distance={18mm}, thick, main/.style = {draw, circle}] 
\node[main,  draw , circle , inner sep =9 pt , thin] (1) {$c_1$}; 
\node[main,  draw , circle , inner sep =9 pt , thin] (2) [above right of=1] {$1$};
\node[main,  draw , circle , inner sep =9 pt , thin] (3) [below right of=2] {$c_2$};
\node[main,  draw , circle , inner sep =9 pt , thin] (4) [above right of=3] {$2$};
\node[main,  draw , circle , inner sep =9 pt , thin] (5) [below right of=4] {$c_3$};
\node[main,  draw , circle , inner sep =9 pt , thin] (6) [above right of=5] {$3$};
\node[main,  draw , circle , inner sep =9 pt , thin] (7) [below right of=6] {$c_4$};
\node[main,  draw , circle , inner sep =9 pt , thin] (8) [above right of=7] {$...$};
\node[main,  draw , circle , inner sep =4 pt , thin] (9) [below right of=8] {$c_{N-1}$};
\node[main,  draw , circle , inner sep =2 pt , thin] (10) [above right of=9] {$N-1$};
\node[main,  draw , circle , inner sep =8 pt , thin] (11) [below right of=10] {$c_{N}$};
\node[main,  draw , circle , inner sep =7 pt , thin] (12) [above right of=11] {$N$};

\draw[->] (1) -- node[sloped, above] {p} (2); 
\draw[->] (2) -- node[sloped, above] {1} (3); 
\draw[->] (1) -- node[sloped, below] {1-p} (3); 
\draw[->] (3) -- node[sloped, above] {p} (4); 
\draw[->] (3) -- node[sloped, below] {1-p} (5); 
\draw[->] (4) -- (5); 
\draw[->] (4) -- node[sloped, above] {1} (5);
\draw[->] (5) -- node[sloped, above] {p} (6); 
\draw[->] (6) -- node[sloped, above] {1} (7); 
\draw[->] (5) -- node[sloped, below] {1-p} (7); 
\draw[dashed,->] (7) -- (8);
\draw[dashed,->] (7) -- (9);
\draw[dashed,->] (8) -- (9);
\draw[->] (9) -- node[sloped, above] {p} (10); 
\draw[->] (9) -- node[sloped, below] {1-p} (11); 
\draw[->] (10) -- node[sloped, above] {1} (11); 
\draw[->] (11) -- node[sloped, above] {1} (12);

\end{tikzpicture}

A trajectory is described as an ordered tuple of probabilities, for instance the tuple  $(p, 1-p, p)$ describes a trajectory which passes thought the first node, avoids the second node (2), and attains the third module, all in two time units.  A similar trajectory with attains the third node in two time units is  $(1-p, p, p)$. The probability of attaining the third node in two time units is therefore $p.(1-p). p + (1-p).p.p = p(1-p).p =  p \begin{pmatrix} 2  \\
    1  \\ \end{pmatrix}  p(1-p)$. The probabitity of having attained a node at position k at time point $m$  is given by:

$$ p(x = k, t = m \tau^+) =  p \begin{pmatrix}
    k-1  \\
    m  \\
\end{pmatrix} p^m  (1-p)^{k-1-m} ~~~~~~~~~~  m \leq k-1  $$

The probability of attaining the terminal node before the time m is:

$$ p(x = N, t \leq m \tau^+) = \sum_{w =0}^m \begin{pmatrix}
    N-1  \\
    w  \\
\end{pmatrix} p^w(1-p)^{N-1-w} $$

The sum of the last two probabilities is equal to one:

$$ \sum_{k = m+1}^{N-1} p(x = k, t = m \tau^+)  +  p(x = N, t \leq m \tau^+)  = 1 $$

We then deduce the formula:

\begin{equation}
\boxed{
  \sum_ {k = m+1}^{N-1}   \begin{pmatrix}
    k-1  \\
    m  \\
\end{pmatrix} p^{m+1}  (1-p)^{k-1-m}  +  \sum_{w =0}^m \begin{pmatrix}
    N-1  \\
    w  \\
\end{pmatrix} p^w(1-p)^{N-1-w}  = 1
}
\end{equation}

\begin{center}
\includegraphics[scale = .35]{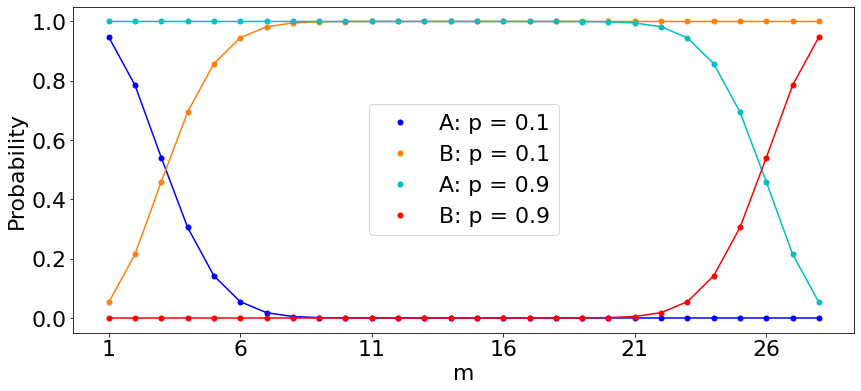}
\captionof{figure}{Probability decomposed as sum of A (first summand) and B (second), as function of time m, for p = 0.1 and p = 0.9. N = 28.}
\end{center}

\section{Two nodes}

\begin{tikzpicture}[node distance={23mm}, thick, main/.style = {draw, circle}] 
\node[main,  draw , circle , inner sep =9 pt , thin] (1) {$c_1$}; 
\node[main,  draw , circle , inner sep =9 pt , thin] (2) [above right of=1] {$1$};
\node[main,  draw , circle , inner sep =9 pt , thin] (3) [below right of=2] {$c_2$};
\node[main,  draw , circle , inner sep =9 pt , thin] (4) [above right of=3] {$2$};
\node[main,  draw , circle , inner sep =9 pt , thin] (5) [below right of=4] {$c_3$};

\node[main,  draw , circle , inner sep =9 pt , thin] (8) [above right of=5] {$...$};
\node[main,  draw , circle , inner sep =4 pt , thin] (9) [below right of=8] {$c_{N-1}$};
\node[main,  draw , circle , inner sep =2 pt , thin] (10) [above right of=9] {$N-1$};
\node[main,  draw , circle , inner sep =8 pt , thin] (11) [below right of=10] {$c_{N}$};
\node[main,  draw , circle , inner sep =7 pt , thin] (12) [above right of=11] {$N$};
\node[main,  draw , circle , inner sep =8 pt , thin] (13) [above of=2]{$1'$}; 
\node[main,  draw , circle , inner sep =9 pt , thin] (14) [above of=4] {$2'$};
\node[main,  draw , circle , inner sep =10 pt , thin] (15) [above of=8] {$...$};
\node[main,  draw , circle , inner sep =2 pt , thin] (16) [above of=10] {$N-1'$};

\draw[->] (1) -- node[sloped, above] {$p_1$} (2); 
\draw[->] (2) -- node[sloped, below] {$1-p_2$} (3);
\draw[->] (1) -- node[sloped, below] {$1-p_1$} (3); 
\draw[->] (3) -- node[sloped, above] {$p_1$} (4); 
\draw[->] (3) -- node[sloped, below] {$1-p_1$} (5); 
\draw[->] (4) -- (5); 
\draw[->] (4) -- node[sloped, below] {$1-p_2$} (5);
\draw[dashed, ->] (5) -- node[sloped, above] {$p_1$} (8); 
\draw[dashed, ->] (5) -- node[sloped, below] {$1-p_1$} (9); 
\draw[dashed,->] (8) -- (9);
\draw[->] (9) -- node[sloped, above] {$p_1$} (10); 
\draw[->] (9) -- node[sloped, below] {$1-p_1$} (11); 
\draw[->] (10) -- node[sloped, below] {$1-p_2$} (11); 
\draw[->] (11) -- node[sloped, above] {$1$} (12);

\draw[->] (2) -- node[left] {$p_2$} (13);
\draw[->] (4) -- node[left] {$p_2$} (14);
\draw[->] (8) -- node[left] {$p_2$} (15);
\draw[->] (10) -- node[left] {$p_2$} (16);

\draw[->] (13) -- node[sloped, above] {$1$} (3);
\draw[->] (14) -- node[sloped, above] {$1$} (5);
\draw[->] (15) -- node[sloped, above] {$1$} (9);
\draw[->] (16) -- node[sloped, above] {$1$} (11);

\end{tikzpicture}

$$ p(x = k, t = m \tau^+) =  \sum\limits_{\substack{ 2a + b +j = m \\ a,b \in  [\![  0,m ]\!]  \\ j \in \{ 0,1\}   }} \begin{pmatrix}
    k-1  \\
    a  \\
\end{pmatrix} \begin{pmatrix}
    k-1-a  \\
    b  \\
\end{pmatrix}   (p_1 p_2^j)(p_1 p_2)^a (p_1(1-p_2))^b  (1-p_1)^{k-a-b-1} $$

x = N: terminal node.

$$ p(x = N, t \leq m \tau^+) =   \sum_{w =0}^m \sum\limits_{\substack{ 2a + b = w \\ a,b \in  [\![  0,m ]\!]}}  \begin{pmatrix}
    N-1  \\
    a  \\
\end{pmatrix} \begin{pmatrix}
    N-1-a  \\
    b  \\
\end{pmatrix} (p_1 p_2)^a (p_1(1-p_2))^b  (1-p_1)^{N-a-b-1} $$

Therefore:
$$ \sum_ {k = 1 + \lfloor \frac{m}{2} \rfloor}^{N-1}  p(x= k, t = m\tau^+)  + p(x = N, t \leq m \tau^+)  = 1$$

\begin{equation}
\begin{split}
  1 =&   \sum_ {k = 1 + \lfloor \frac{m}{2} \rfloor}^{N-1} \sum\limits_{\substack{ 2a + b +j = m \\ a,b \in  [\![  0,m ]\!]  \\ j \in \{ 0,1\}   }} \begin{pmatrix}
    k-1  \\
    a  \\
\end{pmatrix} \begin{pmatrix}
    k-1-a  \\
    b  \\
\end{pmatrix}   (p_1 p_2^j)(p_1 p_2)^a (p_1(1-p_2))^b  (1-p_1)^{k-a-b-1} \\ 
+& \sum_{w =0}^m \sum\limits_{\substack{ 2a + b = w \\ a,b \in  [\![  0,m ]\!]}}  \begin{pmatrix}
    N-1  \\
    a  \\
\end{pmatrix} \begin{pmatrix}
    N-1-a  \\
    b  \\
\end{pmatrix} (p_1 p_2)^a (p_1(1-p_2))^b  (1-p_1)^{N-a-b-1}   
\end{split}
\end{equation}
\section{Three nodes}
\resizebox{1\textwidth}{!}{
\begin{tikzpicture}[node distance={23mm}, thick, main/.style = {draw, circle}]
\node[main,  draw , circle , inner sep =9 pt , thin] (1) {$c_1$}; 
\node[main,  draw , circle , inner sep =9 pt , thin] (2) [above right of=1] {$1$};
\node[main,  draw , circle , inner sep =9 pt , thin] (3) [below right of=2] {$c_2$};
\node[main,  draw , circle , inner sep =9 pt , thin] (4) [above right of=3] {$2$};
\node[main,  draw , circle , inner sep =9 pt , thin] (5) [below right of=4] {$c_3$};

\node[main,  draw , circle , inner sep =9 pt , thin] (8) [above right of=5] {$...$};
\node[main,  draw , circle , inner sep =4 pt , thin] (9) [below right of=8] {$c_{N-1}$};
\node[main,  draw , circle , inner sep =2 pt , thin] (10) [above right of=9] {$N-1$};
\node[main,  draw , circle , inner sep =8 pt , thin] (11) [below right of=10] {$c_{N}$};
\node[main,  draw , circle , inner sep =7 pt , thin] (12) [above right of=11] {$N$};
\node[main,  draw , circle , inner sep =8 pt , thin] (13) [above of=2]{$1'$}; 
\node[main,  draw , circle , inner sep =9 pt , thin] (14) [above of=4] {$2'$};
\node[main,  draw , circle , inner sep =10 pt , thin] (15) [above of=8] {$...$};
\node[main,  draw , circle , inner sep =2 pt , thin] (16) [above of=10] {$N-1'$};

\node[main,  draw , circle , inner sep =8 pt , thin] (17) [above of=13] {$1''$};
\node[main,  draw , circle , inner sep =8 pt , thin] (18) [above of=14] {$2''$};
\node[main,  draw , circle , inner sep =11 pt , thin] (19) [above of=15] {$...$};
\node[main,  draw , circle , inner sep =3 pt , thin] (20) [above of=16] {$N-1''$};

\draw[->] (1) -- node[sloped, above] {$p_1$} (2); 
\draw[->] (2) -- node[sloped, below] {$1-p_2$} (3);
\draw[->] (1) -- node[sloped, below] {$1-p_1$} (3); 
\draw[->] (3) -- node[sloped, above] {$p_1$} (4); 
\draw[->] (3) -- node[sloped, below] {$1-p_1$} (5); 
\draw[->] (4) -- (5); 
\draw[->] (4) -- node[sloped, below] {$1-p_2$} (5);
\draw[dashed, ->] (5) -- node[sloped, above] {$p_1$} (8); 
\draw[dashed, ->] (5) -- node[sloped, below] {$1-p_1$} (9); 
\draw[dashed,->] (8) -- (9);
\draw[->] (9) -- node[sloped, above] {$p_1$} (10); 
\draw[->] (9) -- node[sloped, below] {$1-p_1$} (11); 
\draw[->] (10) -- node[sloped, below] {$1-p_2$} (11); 
\draw[->] (11) -- node[sloped, above] {$1$} (12);

\draw[->] (2) -- node[left] {$p_2$} (13);
\draw[->] (4) -- node[left] {$p_2$} (14);
\draw[->] (8) -- node[left] {$p_2$} (15);
\draw[->] (10) -- node[left] {$p_2$} (16);

\draw[->] (13) -- node[pos=0.25, sloped, below] {$1-p_3$} (3);
\draw[->] (14) -- node[pos=0.25, sloped, below] {$1-p_3$} (5);
\draw[->] (15) -- node[pos=0.25, sloped, below] {$1-p_3$} (9);
\draw[->] (16) -- node[pos=0.25, sloped, below] {$1-p_3$} (11);

\draw[->] (13) -- node[left] {$p_3$} (17);
\draw[->] (14) -- node[left] {$p_3$} (18);
\draw[->] (15) -- node[left] {$p_3$} (19);
\draw[->] (16) -- node[left] {$p_3$} (20);

\draw[->] (17) -- node[sloped, above] {$1$} (3);
\draw[->] (18) -- node[sloped, above] {$1$} (5);
\draw[->] (19) -- node[sloped, above] {$1$} (9);
\draw[->] (20) -- node[sloped, above] {$1$} (11);

\end{tikzpicture}
}
$$ $$

We give the rule for 3 nodes, which  the suggests the form of the general formula, for arbitrary number of nodes.

\begin{align*} p(x = k, t = m \tau^+) =  \sum\limits_{\substack{ 3a + 2b + c +j  + jl = m \\ a,b,c \in  [\![  0,m ]\!]  \\ j,l \in \{ 0,1\}   }} &\begin{pmatrix}
    k-1  \\
    a  \\
\end{pmatrix} \begin{pmatrix}
    k-1-a  \\
    b  \\
\end{pmatrix} \begin{pmatrix}
    k-1-a-b  \\
    c  \\
\end{pmatrix}   (p_1 p_2^j p_3^{jl})(p_1 p_2 p_3)^a (p_1 p_2(1-p_3))^b... \\ &~~~~~~...(p_1 (1-p_2))^c    (1-p_1)^{k-a-b-c-1} 
\end{align*}

\begin{align*} p(x = k, t \leq m \tau^+) = \sum_{w =0}^m  \sum\limits_{\substack{ 3a + 2b + c   = m \\ a,b,c \in  [\![  0,m ]\!      ]}} &\begin{pmatrix}
    k-1  \\
    a  \\
\end{pmatrix} \begin{pmatrix}
    k-1-a  \\
    b  \\
\end{pmatrix} \begin{pmatrix}
    k-1-a-b  \\
    c  \\
\end{pmatrix}   (p_1 p_2 p_3)^a (p_1 p_2(1-p_3))^b... \\ &~~~~~~...(p_1 (1-p_2))^c    (1-p_1)^{k-a-b-c-1} \\
\end{align*}

\section{Random walk in two dimensions}
  Consider a point initially at the origin, doing random jumps on $\mathbb{Z}^2$, on the possible directions, with probabilities $p_{1,2,3,4}: p_1 + p_2 + p_3 + p_4 = 1 $.
   The jumps can be of arbitrary size and direction. Call $T_i$ the travelling time ($T_i \in  \mathbb{N}^* $). Considering a certain time point $m$ ($m \in  \mathbb{N}^* $), and a pair of coordinates $(h,v) \in \mathbb{Z}^2$. 
Let's study the case: 
    $p_1:  (\rightarrow, J_1, T_1)~~ p_2:  (\leftarrow, J_2, T_2) ~~ p_3:  (\uparrow, J_3, T_3) ~~
  p_4: (\downarrow,  J_4, T_4) $.

Then: 

$$1 =  \sum\limits_{\substack{ h \in [\![J_2.m, J_1.m]\!] \\ v \in [\![J_4.m, J_3.m]\!] }}  p \Big( \begin{pmatrix}
    x  \\
    y  \\
\end{pmatrix} = \begin{pmatrix}
    h  \\
    v  \\
\end{pmatrix},~t = m \Big) $$

\begin{equation}
 1 =   \sum\limits_{\substack{ h \in [\![J_2.m, J_1.m]\!] \\ v \in [\![J_4.m, J_3.m]\!] }} 
  \sum\limits_{\substack{ J_1.a + J_2.b = h     \\ J_3.c + J_4.d = v \\ a.T_1+b.T_2+c.T_3+d.T_4 = m  }} \begin{pmatrix}
    a+b  \\
    a  \\
\end{pmatrix} \begin{pmatrix}
    c+d  \\
    c  \\
\end{pmatrix} p_1^a . p_2^b . p_3^c.  p_4^d
\end{equation}

In the second summand, we iterate over  all pair of values $(a, b, c, d)$ which are solutions in $\mathbb{N}^4$, of this   equations system (two spatial equations and one temporal equation).

\subsection{Simple random walk in 1D}

   Starting from the origin, consider a point which moves $+1$ space unit with probability $p$ and $-1$ with probability $1-p$. The jumps last both one time unit. After $m$ time steps, the point is situated in the interval $[-m,m]$. 
We will only observe the point position at "even" moments in time ($m = 2n, n \in \mathbb{N}$ even), and therefore we only consider even positions ($k =2l, l \in \mathbb{Z}$ even). The probability of being in odd positions is nil at even times.  We find the combinatorial sum:
\begin{equation}
 1 = \sum\limits_{\substack{ k \in [\![-m, m]\!] \\ k ~even      }}  p(x = k, t = m)  \iff 1 = \sum\limits_{\substack{ k \in [\![-m, m]\!] \\ k ~even      }}   \begin{pmatrix}
         m  \\
         \frac{m+k}{2}  \\
         \end{pmatrix} p^{ \frac{m+k}{2}}(1-p)^{ \frac{m-k}{2}} 
\end{equation}
\subsection{Simple chain with delay and a terminal node}

We can consider in 1D a chain in which elements go wright by one spatial unit  with probability $p$ (with a waiting time of 1 time unit), and go wright, with jumps of one spatial unit,  with a probability of $1-p$, instantaneously.  We find back the formula of the chain presented in the second paragraph.

\subsection{Simple 2D chain with delay and a terminal barrier}

 We consider that a point has the following dynamics: 
 
 $p_1 (\rightarrow, J_1 = 1, T_1 = 1)$, 
 $~p_2 (\rightarrow, J_2 = 1, T_2 = 0)$,
  $~p_3 (\uparrow, J_3 = 1, T_3 = 1)$, 
 $~p_4 (\uparrow, J_4 = 1, T_4 = 0)$.

$$p \Big( \begin{pmatrix}
    x  \\
    y  \\
\end{pmatrix} = \begin{pmatrix}
    h  \\
    v  \\
\end{pmatrix},~t = m \Big)  = \sum_{\substack{a+c = m \\ a \in [\![0, h]\!] \\ c \in [\![0, v]\!] }}      \begin{pmatrix}
    h  \\
    a  \\
\end{pmatrix}\begin{pmatrix}
    v  \\
    c  \\
\end{pmatrix} p_1^a . p_2^{h-a} . p_3^c.  p_4^{v-c}   $$ 
 
 Lets consider that at there is a barrier $\mathcal{B}$ delimited by the coordinates:
 
 $$\mathcal{B} = \{  \begin{pmatrix}
    N  \\
    v  \\
\end{pmatrix},v \in [\![0,N-1 ]\!] ~or  
\begin{pmatrix}
    h  \\
    N  \\
\end{pmatrix},h \in [\![0,N-1]\!]
  \}  $$ 
 
Therefore, at time m, the probability of being in a certain position (that has not reached the barrier), plus the probability of having reached the barrier is equal to one.

\begin{equation}
 \sum_{ \substack{    \begin{pmatrix}
    h  \\
    v  \\
\end{pmatrix} \not\in \mathcal{B} \\ m\leq h+v     } }    \sum_{\substack{a+c = m \\ a \in [\![0, h]\!] \\ c \in [\![0, v]\!] }}      \begin{pmatrix}
    h  \\
    a  \\
\end{pmatrix}\begin{pmatrix}
    v  \\
    c  \\
\end{pmatrix} p_1^a . p_2^{h-a} . p_3^c.  p_4^{v-c}   +
  \sum_{\substack{\begin{pmatrix}
    h  \\
    v  \\
\end{pmatrix} \in \mathcal{B}\\ m\leq h+v     }}    \sum_{\substack{a+c = 0 \\ a \in [\![0, h]\!] \\ c \in [\![0, v]\!] }}^{m}      \begin{pmatrix}
    h  \\
    a  \\
\end{pmatrix}\begin{pmatrix}
    v  \\
    c  \\
\end{pmatrix} p_1^a . p_2^{h-a} . p_3^c.  p_4^{v-c}  = 1  
\end{equation}

\begin{center}
\newcommand*{\xMin}{0}%
\newcommand*{\xMax}{6}%
\newcommand*{\yMin}{0}%
\newcommand*{\yMax}{6}%

\begin{tikzpicture}[scale=0.50]
    \foreach \i in {\xMin,...,\xMax} {
        \draw [very thin,gray] (\i,\yMin) -- (\i,\yMax)  node [below] at (\i,\yMin) {$\i$};
    }
    \foreach \i in {\yMin,...,\yMax} {
        \draw [very thin,gray] (\xMin,\i) -- (\xMax,\i) node [left] at (\xMin,\i) {$\i$};
    }

\draw [very thick, cyan, step=1.0cm,xshift=0cm, yshift=0cm] (0,0) grid+(5,5); 

\draw [very thick, blue] (5,0) -- (5,5);
\draw [very thick, blue] (0,5) -- (5,5);

\draw[thick,->] (0,0) -- (7,0) node[anchor=north west] {x};
\draw[thick,->] (0,0) -- (0,7) node[anchor=south east] {y};

\end{tikzpicture}
\end{center}

\section{Conclusion}

The formula deduced in equation (1) has  some ressemblance with Gosper's formula  \citep{beeler1972hakmem}.
$$1 = \sum_{k = 0}^{N-1} \begin{pmatrix}
    N+k-1  \\
    k  \\
\end{pmatrix} \Big(p^{N}(1-p)^k  + (1-p)^{N}p^k \Big) $$

Consider a coin flip: it has a probability p of falling heads, and a probability 1-p of being tails. This formula reflects that if the flip the coin until one of its sides comes up N times, the sum of the probability all the probable states $\mathcal{S}$ is equal to one.
$$\mathcal{S} = \{(N,0), (N, 1), ..., (N,N-1) \} \bigcup \{(0,N), (1,N), ..., (N-1, N) \} $$

  In equation (1), however, we "take the picture", not after N heads or tails like in the Gosper's formula, but at an arbitraty time m, ($m \leq N)$.
After m flips we compute the probabilities of having reached different states. In fact the sum of the exponents in Gospers Formula is equal to $N+k, ~~k\in [0,N-1]$, whereas in (1) it is equal to k ($k \in [\![m+1, N-1 ]\!] $).

\bibliographystyle{unsrtnat}
\bibliography{references}

\end{document}